\documentclass[a4paper,12pt]{amsart}
\usepackage[utf8]{inputenc}

\usepackage{amsthm}
\usepackage{verbatim}
\usepackage{amsmath}
\usepackage{amssymb}
\usepackage{enumerate}

\allowdisplaybreaks
\numberwithin{equation}{section}
\numberwithin{figure}{section}

\theoremstyle{plain}
\newtheorem{theorem}{Theorem}[section]
\newtheorem{proposition}[theorem]{Proposition}
\newtheorem{lemma}[theorem]{Lemma}
\newtheorem{corollary}[theorem]{Corollary}
\newtheorem{example}[theorem]{Example}
\newtheorem{remark}[theorem]{Remark}
\newtheorem{definition}[theorem]{Definition}
\newtheorem{conjecture}{Conjecture}

\def\BET{\begin{theorem}}
\def\ENT{\end{theorem}}
\def\BEP{\begin{proposition}}
\def\ENP{\end{proposition}}
\def\BEL{\begin{lemma}}
\def\ENL{\end{lemma}}
\def\BEC{\begin{corollary}}
\def\ENC{\end{corollary}}
\def\BEE{\begin{example}\rm}
\def\ENE{\end{example}}
\def\BER{\begin{remark} \rm}
\def\ENR{\end{remark}}
\def\BED{\begin{definition} \rm}
\def\END{\end{definition}}
\def\BECJ{\begin{conjecture}}
\def\ENCJ{\end{conjecture}}

\def\bea{\begin{eqnarray}}
\def\eea{\end{eqnarray}}
\def\bean{\begin{eqnarray*}}
\def\eean{\end{eqnarray*}}
\def\beq{\begin{equation}}
\def\eeq{\end{equation}}
\def\beal{\begin{align*}}
\def\eeal{ \end{align*} }

\def\bbC{{\mathbb C}}
\def\bbD{{\mathbb D}}

\def\bbN{{\mathbb N}}

\begin{document}


\title[Vanishing Bergman kernels]{Vanishing Bergman kernels on the disk}



\author{Antti Per\"al\"a}
\address{Department of Physics and Mathematics, University of Eastern Finland, P.O.Box 111, FI-80101 Joensuu, Finland}
\email{antti.perala@uef.fi}

\thanks{This research is supported by the Academy of Finland project no. 268009.}

\keywords{Bergman space, Bergman kernel}
\subjclass[2010]{30H20, 32A36}
\begin{abstract}
We discuss topics related to zeroes of the Bergman kernels, and present a method for generating Bergman kernels with arbitrarily, but finitely, many zeroes. It is also shown that a Bergman kernel induced by a radial weight on the unit disk cannot have infinitely many zeroes. Similar questions for the Segal-Bargmann spaces of the complex plane are briefly discussed.
\end{abstract}

\maketitle

\section{Introduction.}
Denote by $\bbD$ the open unit disk of the complex plane. An integrable function $\omega:\bbD \to [0,\infty)$ is called a weight. We are mainly interested in radial weights, that is, the weights $\omega$ satisfying $\omega(z)=\omega(|z|)$ for all $z \in \bbD$. Denote by $dA(z)=\pi^{-1}dxdy$ ($z=x+iy$) the normalized Lebesgue measure on the disk, and define the Bergman space $A_\omega^p$ ($0<p<\infty$) to consist of analytic functions $f:\bbD\to\bbC$ with
\begin{align*}
\|f\|_{A_\omega^p}&=\left(\int_\bbD |f(z)|^p\omega(z)dA(z)\right)^{1/p}\\
&=\left(\pi^{-1}\int_0^1 \int_0^{2\pi}|f(re^{i\theta})|^p d\theta \omega(r)rdr\right)^{1/p}<\infty.
\end{align*}
It is customary, even for $0<p<1$, to talk about norms of elements of $A_\omega^p$, as well as norms of their operators. The reader should not be confused by this slight abuse of terminology.
Of course, when $p=2$, the norm above is given in form of the inner product $\langle f,f\rangle_\omega^{1/2}$, where
$$\langle f,g\rangle_\omega=\int_\bbD f(z)\overline{g(z)}\omega(z)dA(z).$$
A weight $\omega$ is called a Bergman weight, if for every $0<p<\infty$ convergence in $A^p_\omega$ implies uniform convergence in compact sets. From this it follows that for every $z \in \bbD$, the functional
$$f\mapsto f(z)$$
is bounded on $A_\omega^p$. By the Riesz representation theorem, there then exists a unique $B_z^\omega \in A_\omega^2$ such that
$$f(z)=\langle f|B_z^\omega \rangle_\omega$$ for $f \in A_\omega^2$. The function $B_z^\omega$ is called the Bergman kernel.

We will assume throughout that $\omega$ is a Bergman weight. We will also employ the following convention. The weights $\omega$ and $\nu$ are equivalent, if $A^2_\omega=A^2_\nu$ as sets. In this case, we write $\omega \sim \nu$. The condition $\omega \sim \nu$ is weaker than requiring that the weights are pointwise comparable; it is very well possible that $\omega$ has an $L^1$ singularity, while $\nu$ vanishes near the origin, and still $\omega \sim \nu$.

Note that if $\omega$ is radial, then for every $z \in \bbD$ the function $B_z^\omega$ has the series representation
\begin{equation}\label{srs}
B_z^\omega(\xi)=\sum_{n=0}^\infty \frac{1}{\omega_n}(\overline{z}\xi)^n,
\end{equation}
where the numbers $\omega_n=2\int_0^1 r^{2n+1}\omega(r)dr$ are the moments of $\omega$. In this case $B_z^\omega$ is actually analytic in a larger disk than $\bbD$, thus the reproducing formula makes sense if $f\in A^1_\omega$. In fact, by using a limit argument, the formula can be justified for any analytic $f$. A cornerstone of the theory of reproducing kernels is the work of N. Aronszajn \cite{Aro}.

In this paper we are interested on the zeroes of $B_z^\omega$, in particular, whether they exist. The vanishing of Bergman kernels on domains of $\bbC^n$ is known as the Lu Qi-Keng's problem \cite{LQK}. It has been studied extensively by H. P. Boas, S. Fu, E. J. Straube, M. Engli\u s and others; we list here the following papers: \cite{Boa}, \cite{BFS}, \cite{Boa2} and \cite{Eng}. We also mention the paper of Y. E. Zeytuncu \cite{Zey}, where a weight $\lambda$ on $\bbD$ is constructed so that $\lambda \sim 1$ and $B_z^\lambda$ has a zero.

The paper is composed as follows. We continue the discussion about non-vanishing kernels in the next section, with emphasis on motivating the topic. The third section contains the three main results of the paper. First we show that for any $\alpha>1$, there exists a radial weight $\omega_\alpha$, which is equivalent to the standard weight $\nu_\alpha$ (see the next section), such that the Bergman kernel $B^{\omega_\alpha}_z$ has a zero for some $z \in \bbD$. The second main result shows that it is possible to find a radial weight, which is equivalent to some standard weight, such that the associated kernel has any finite prescribed number of zeroes. Finally, we show that for a radial weight, the Bergman kernel cannot have infinite number of zeroes. In the last section of the paper, we briefly comment of analogues of the main results for the Segal-Bargmann spaces of the complex plane.

As is standard, by $a \lesssim b$ we mean that, under appropriate constraints, there exists $C=C(\cdot)>0$ such that $a \leq Cb$. The relation $a\gtrsim b$ is defined analogously, and if both $a\lesssim b$ and $a \gtrsim b$ hold, we write $a\asymp b$.

\section{Motivating remarks.}

A central family of weights is given by $\nu_\alpha(z)=(\alpha+1)(1-|z|^2)^{\alpha}$, where $\alpha>-1$. These are the so-called standard weights. They are studied in numerous papers and covered in many textbooks, such as \cite{HKZ} and \cite{Z}. It is well-known that
$$B_z^{\nu_\alpha}(\xi)=\frac{1}{(1-\overline{z}\xi)^{2+\alpha}}.$$ The spaces $A_{\nu_\alpha}^p$ have the nice property that
$$f\mapsto (f\circ \varphi_z)\varphi_z'^{(2+\alpha)/p}$$ is an isometric isomorphism $A_{\nu_\alpha}^p\to A_{\nu_\alpha}^p$. Here the map
$$\varphi_z(\xi)=\frac{z-\xi}{1-\overline{z}\xi}$$ is the M\"obius automorphism of the disk, which interchanges $0$ and $z$. In \cite{Vuk}, D. Vukoti\'c used this fact to obtain the sharp bound
\begin{equation}\label{vuk}
|f(z)|\leq (1-|z|^2)^{-(2+\alpha)/p}\|f\|_{A^p_{\nu_{\alpha}}},
\end{equation}
which holds for every point $z \in \bbD$.

\BER We note here that the case $\alpha=-1$ for $\nu_\alpha$ can often (but not always, see \cite{PRMem}, for instance) be interpreted as a result about the Hardy spaces $H^p$. In fact, setting $\alpha=-1$ in the inequality \eqref{vuk}, one recovers the corresponding sharp bound for $H^p$. Of course, since $|f|^p$ is subharmonic, the weight $(1-|z|^2)^{-1}$ would induce a trivial Bergman space. However, consider the following construction, which we present for the reader's amusement. Let $\eta$ be a radial, continuous weight that is a priori only locally integrable. To an analytic function, we assign, for every $R \in (0,1)$, the quantity
$$n_\eta^p(f,R):=\frac{\int_{\{|z|<R\}} |f(z)|^p \eta(z)dA(z)}{\int_{\{|z|<R\}}\eta(z)dA(z)}.$$
When $\eta \in L^1$, we see that $n_\eta^p(f,R)\to \|f\|_{A^p_{\widetilde{\eta}}}^p$ as $R\to 1^-$, where $\widetilde{\eta}$ is just $\eta$ normalized to a probability measure. Even if $\eta \in L^1_{loc}\setminus L^1$, it is clear that
\begin{equation}\label{limsup}
\limsup_{R\to 1^-}n_\eta^p(f,R)<\infty
\end{equation} for bounded analytic functions $f$, so it makes sense to consider the set of analytic functions for which \eqref{limsup} holds. Note that if $f$ is not the zero function, by subharmonicity of $|f|^p$, we see that both the numerator and the denominator of $n_\eta^p(f,R)$ tend to infinity as $R\to 1^-$. Now, writing in polar coordinates and applying L'H\^opital's rule to get rid of the integral, we see
$$n_\eta^p(f,R)=\frac{\pi^{-1}\int_0^R\int_0^{2\pi}|f(re^{i\theta})|^p d\theta r\eta(r)}{2\int_0^R r\eta(r)dr}\to \lim_{R\to 1^-} \frac{1}{2\pi}\int_0^{2\pi}|f(Re^{i\theta})|^pd\theta,$$
which is finite if and only if $f \in H^p$.

In view of this remark, the Bergman spaces induced by any $(1-|z|^2)^{\alpha}$ for $\alpha\leq -1$ would coincide with the Hardy space of same integrability exponent. From another viewpoint, it is sometimes more convenient to think of the case $\alpha=-2$ as the analytic Besov space.
\ENR

For the purpose of the present paper, yet another nice property of the weights $\nu_\alpha$ is that the Bergman kernels $B_z^{\nu_\alpha}$ have no zeroes. This allows arbitrary powers of the kernels while preserving analyticity, which is tremendously useful in operator theory. The upper bound in the following result is implicit in the paper \cite{Z2000}, but we present it here for the convenience of the reader and show that it is sharp.

\BET\label{sharp}
Let $\omega$ be a (not necessarily radial) weight, $1\leq p<\infty$, and $z \in \bbD$ a point for which $B_z^\omega$ has no zeroes. Then the norm of the functional
\begin{equation}
\lambda_z(f)=f(z)
\end{equation} equals $B_z^\omega(z)^{1/p}$.
\ENT

\begin{proof}
Let $q$ be the conjugate exponent of $p$; $1/p+1/q=1$. Consider the function
$$F_z^{\omega,q}(\xi)=B_z^\omega(\xi)\overline{B_z^\omega(\xi)}^{2/q-1}B_z^\omega(z)^{1-2/q},$$
which is well-defined since $B_z^\omega$ has no zeroes.
Note that
$$|F_z^{\omega,q}(\xi)|^q = |B_z^\omega(\xi)|^2 B_z^\omega(z)^{q-2}$$ so that
If $f \in A^p_\omega$, then
$$\langle f,F_z^{\omega,q}\rangle_\omega=B_z^\omega(z)^{1-2/q}\int_\bbD \left[f(\xi)B_z^\omega(\xi)^{2/q-1}\right]\overline{B_z^\omega(\xi)}\omega(\xi)dA(\xi)=f(z).$$
By H\"older's inequality, we now get
$$|f(z)|\leq B_z^\omega(z)^{1/p}\|f\|_{A^p_\omega}.$$
On the other hand, the function $G_z^{\omega,p}(\xi)=B_z^{\omega}(\xi)^{2/p}$ is analytic and satisfies
$\|G_z^{\omega,p}\|_{A^p_\omega}=B_z^\omega(z)^{1/p}$. Since $G_z^{\omega,p}(z)=B_z^{\omega}(z)^{2/p}$, we see that
$$|G_z^{\omega,p}(z)|= B_z^\omega(z)^{1/p}\|G_z^{\omega,p}\|_{A_\omega^p},$$
so the estimate is sharp.
\end{proof}

\BER Let $1\leq p<\infty$. Since the point-evaluation $\lambda_z$ belongs to the dual of $A^p_\omega$, it can be extended to an element in the dual of $L^p_\omega$ by using the Hahn-Banach theorem. Further, it is well-known that this extension can be expressed as an $L^q_\omega$ function, whose norm agrees with the original dual norm of $\lambda_z$. We remark that the function $F_z^{\omega,q}$ is such an element. In many cases, when $1<p<\infty$, the dual of $A^p_\omega$ is isomorphic to $A^q_\omega$. We mention that the case of regular weights is treated in \cite{PR} Corollary 7, whereas the case of exponential weights (on the disk or the plane) requires that the weight function changes with $p$, see \cite{CP} Corollary 8 and \cite{JPR} Corollary 3.4. As for duality considerations for $0<p\leq 1$, there are the recent papers \cite{PRinfty} and \cite{APR}. The isomorphism involved in these dualities is usually not an isometry, and indeed, as an evidence of this fact, the function $F_z^{\omega,q}$ is very rarely analytic.

It is clear that the proof can very easily be adapted to hold for Bergman spaces on simply connected domains. We also remark that when $0<p<1$, the functions $B_z^\omega(\xi)^{2/p}$ can be used to obtain the same lower bound as in the theorem above, whereas we do not know if the upper bound still holds. Of course, when $p=2$, the above theorem holds even if the kernel functions have zeroes.

This theorem can also be used to obtain the estimate $(1.6)$ of Theorem 1.2 of \cite{GK}, which deals with a similar question for the Segal-Bargmann space. Note that the dependence of $p$ is suppressed by the way the Segal-Bargmann norm depends on $p$. We will address this topic in more detail in the final section of the paper.
\ENR

One downside of the argument in Theorem \ref{sharp} is that it does not give upper bound for the small exponents $0<p<1$. Another (perhaps more serious) shortcoming is that the non-vanishing of $B_z^\omega$ is essential. Fortunately, there are many cases where such vanishing is impossible, the most common being the case of standard weights as mentioned earlier. Due to a result of H. Hedenmalm \cite{DKS} (see also \cite{Hed}) the vanishing of the kernel can be ruled out for (not necessarily radial) weights $\omega$ for which $\log \omega$ is subharmonic, and the zeroes of $\omega$ (if any) are isolated. In fact, in \cite{Z2000} it is shown that the same upper estimate can be obtained for the full range $0<p<\infty$ in this case. As remarked above, this estimate is sharp.

\section{Vanishing kernels for the Bergman spaces.}

In this section, we present a method for generating Bergman kernels with zeroes. Given a radial weight $\omega$, let
$$\omega^*(z)=\int_{|z|}^1\omega(s)\log\left(\frac{s}{|z|}\right)s ds$$ be its associated weight, which is defined on $\bbD\setminus \{0\}$, see \cite{PRMem, PR}. Our first main result is the following.

\BET\label{example}
Let $\omega$ be a radial weight. Then
\begin{equation}\label{star}
B_z^{\omega^*}(\xi)=4\partial_{\xi}\partial_{\overline{z}}B_z^{\omega}(\xi).
\end{equation}
Moreover, if $\alpha>-1$, and $\nu_\alpha(z)=(\alpha+1)(1-|z|^2)^\alpha$ is the standard weight, then the Bergman kernel generated by its associated weight $\nu_\alpha^*$ has zeroes.
\ENT

\begin{proof}
The argument is based on the Littlewood-Paley identity, which itself is a consequence of the Green's formula
$$\langle f,g\rangle_\omega=4\langle f',g'\rangle_{\omega^*}+\omega(\bbD)f(0)\overline{g(0)}.$$

So let $f$ be an analytic function and $F$ its primitive, which vanishes at the origin. The reproducing formula for $B_z^\omega$ gives
$$f(z)=\int_\bbD F(\xi)\overline{\partial_{\overline{z}}B_z^\omega(\xi)}\omega(\xi)dA(\xi).$$
Using the Littlewood-Paley identity, we get
$$f(z)=4\int_\bbD f(\xi)\overline{\partial_\xi\partial_{\overline{z}}B_z^\omega(\xi)}\omega^*(\xi)dA(\xi),$$
where the additional zero-evaluation term vanishes because $F(0)=0$. Since the reproducing kernels are unique, we conclude that the claimed identity for $B_z^{\omega^*}$ holds.

Now, apply this formula for $\omega=\nu_\alpha$ to obtain
$$B_z^{\nu_\alpha^*}(\xi)=4\partial_{\xi}\partial_{\overline{z}}B_z^{\nu_\alpha}(\xi)=4(2+\alpha)\frac{1+(2+\alpha)\overline{z}\xi}{(1-\overline{z}\xi)^{4+\alpha}}.$$
This function vanishes whenever $\overline{z}\xi=1/(2+\alpha)$, which can be achieved for every $\alpha>-1$ by taking $z$ sufficiently close to the boundary of the unit disk.
\end{proof}

\BER
Note that $\nu_\alpha^*\sim \nu_{\alpha+2}$, so the space $A_{\nu_\alpha^*}^2$ can be equipped with an equivalent inner product, whose reproducing kernel does not vanish. We do not know, whether such claim is true about all weights with vanishing kernels. Neither do we know, if the converse of this claim is true, that is, whether every weighted Bergman space carries an equivalent inner product with a vanishing Bergman kernel.

The weights $\omega^*$ always have a mild singularity at the origin. While this does not cause problems with the Bergman space theory, it is worth noting that the weight given as an example in \cite{Zey} is also very much concentrated inside the unit disk. We are currently unable to obtain any conclusive result based on these remarks.
\ENR

Building on Theorem \ref{example}, we can now prove our second main theorem.

\BET\label{prescribe}
For every $n \in \bbN$, there exists a radial weight $\omega^n$, which is equivalent to a standard weight $\nu_\beta$ for some $\beta>-1$, and $z \in \bbD$, such that the Bergman kernel $B^{\omega_n}_z$ has exactly $n$ zeroes in $\bbD$, counting multiplicities.
\ENT

\begin{proof}
The argument is based on iterating the map $\omega \mapsto \omega^*$, so let us define $\omega^{*0}=\omega$ and $\omega^{*n}=(\omega^{*(n-1)})^*$ for $n \in \bbN$. For a radial $\omega$ we write $B^\omega(\zeta)=B_z^\omega(\xi)$, where $\zeta=\overline{z}\xi$. By a straightforward calculation, using \eqref{star}, we see that then
\begin{align*}
B^{\omega^*}_z(\xi)&=4\partial_\xi \partial_{\overline{z}} \sum_{n=0}^\infty \frac{1}{\omega_n}(\overline{z}\xi)^n \\
&=4\sum_{n=1}^\infty \frac{n^2}{\omega_n}(\overline{z}\xi)^{n-1},
\end{align*}
which, in terms of the newly introduced variable $\zeta$, coincides with the identity
\begin{equation}\label{next}
B^{\omega^*}(\zeta)=4(B^\omega)'(\zeta)+4\zeta (B^\omega)''(\zeta).
\end{equation}
We now specialize to the case of standard weights $\nu_\alpha$, where $\alpha>-1$. We already know that the claim holds when $n=1$, and that $B^{\nu_\alpha^*}(\zeta)$ vanishes at precisely one point. Let us write $B_{\alpha,n}(\zeta)=B^{\nu_\alpha^{*n}}(\zeta)$. We will show that
$$B_{\alpha,n}(\zeta)=\frac{p_{\alpha,n}(\zeta)}{(1-\zeta)^{2+\alpha+2n}},$$
where
$$p_{\alpha,n}(\zeta)=\sum_{k=0}^n c_{\alpha,n,k}\zeta^k$$ is a polynomial of degree $n$. Moreover, when $n\geq 1$ is fixed, we have
\begin{equation}\label{asymp}
|c_{\alpha,n,n}|\asymp \alpha^{2n} \text{ and } |c_{\alpha,n,k}| \lesssim \alpha^{2n-1},
\end{equation}
when $\alpha\to \infty$ and $k<n$. The implicit constants are allowed to depend on the parameter $n$.

The claimed property is proved by induction. First, $$p_{\alpha,1}(\zeta)=4(2+\alpha)+4(2+\alpha)^2\zeta,$$ which is clearly of the form we want. Next, suppose that $p_{\alpha,m}$ satisfies \eqref{asymp} for some $m\geq 1$. By using \eqref{next}, and manipulating each term so that they have the same denominator, we have
\begin{align*}
p_{\alpha,m+1}(\zeta)&=4p_{\alpha,m}'(\zeta)(1-\zeta)^2 \\
&+4(2+\alpha+2m)p_{\alpha,m}(\zeta)(1-\zeta) \\
&+4\zeta p_{\alpha,m}''(\zeta)(1-\zeta)^2 \\
&+8(2+\alpha+2m)\zeta p_{\alpha,m}'(\zeta)(1-\zeta)\\
&+4(2+\alpha+2m)(3+\alpha+2m)\zeta p_{\alpha,m}(\zeta).
\end{align*}
From this we readily verify that $p_{\alpha,m+1}$ is a polynomial of degree $m+1$ and we can collect the highest order term as follows:
\begin{align*}
c_{\alpha,m+1,m+1}&=4c_{\alpha,m,m}m\\
&-4c_{\alpha,m,m} (2+\alpha+2m)\\
&+4c_{\alpha,m,m} m(m-1) \\
&-8c_{\alpha,m,m} (2+\alpha+2m) m \\
&+4c_{\alpha,m,m} (2+\alpha+2m)(3+\alpha+2m),
\end{align*}
where the last term has the desired asymptotics in terms of $\alpha$. Thus $|c_{\alpha,m+1,m+1}|\asymp \alpha^{2m+2}$ as $\alpha\to \infty$, which follows from $|c_{\alpha,m,m}|\asymp \alpha^{2m}$ as $\alpha\to \infty$.

As for the lower order terms, note that the asymptotic increment of $\alpha^2$ can only come from the term
$$4(2+\alpha+2m)(3+\alpha+2m)\zeta p_{\alpha,m}(\zeta)$$ in $p_{\alpha,m+1}$. Since the coefficients $c_{\alpha,m,k}$ of $p_{\alpha,m}$ have the asymptotics $|c_{\alpha,m,k}|\lesssim \alpha^{2m-1}$, when $k<m$, whereas the coefficient $c_{\alpha,m,m}$ is already reserved for the highest order term, we can repeat essentially the same argument as before to obtain $|c_{\alpha,m+1,k}|\lesssim \alpha^{2m+1}$ for $k<m+1$. By induction, the estimates in \eqref{asymp} are proven.

Finally, note that the monomial $c_{\alpha,n,n}\zeta^{n}$ has $n$ zeroes (counting multiplicities) in the unit disk $\bbD$. By the asymptotic formula \eqref{asymp}, if $\alpha$ is large enough, we have
$$|c_{\alpha,n,n}\zeta^n|>\left|\sum_{k=0}^{n-1}c_{\alpha,n,k}\zeta^k \right|,$$
when $|\zeta|=1$. By applying Rouche's theorem, we deduce that such $p_{\alpha,n}$ has exactly $n$ zeroes (counting multiplicities) inside $\bbD$. Finally, since
$$B_{\alpha,n}(\zeta)=\frac{p_{\alpha,n}(\zeta)}{(1-\zeta)^{2+\alpha+2n}},$$
we have the same property for $B_{\alpha,n}$.

Let now $\zeta_n$ be a zero of $B_{\alpha,n}$ with the highest modulus. If $|\zeta_n|<|z|<1$, then it is easy to see that $B_z^{\nu_\alpha^{*n}}$ has exactly $n$ zeroes. So, we put $\omega^n=\nu_\alpha^{*n}$. To complete the proof, note that $\nu_\alpha^{*n}\sim \nu_{\alpha+2n}$, which follows, for instance by iterating Lemma 1.7 of \cite{PRMem} or part $(iv)$ of Lemma A in \cite{PR}.
\end{proof}

\BER
Note that the theorem does not hold for all $\alpha>-1$ and for all $n$, but control on the size of $\alpha$ is needed. For instance, $p_{0,2}(\zeta)$ (see the preceding proof) is a constant multiple of $3\zeta^2+6\zeta+1$, whose roots are $\zeta_1=-1+\sqrt{24}/6$ and $\zeta_2=-1-\sqrt{24}/6$, out of which only $\zeta_1$ belongs to the unit disk.
\ENR

We have so far shown that a Bergman kernel induced by a radial weight can have any finite number of zeroes, and in fact the radial weight in question can be chosen to be equivalent to a standard weight. The following theorem shows that the possibility of having infinite number of zeroes can be excluded for all radial weights.

\BET\label{inf}
Let $\omega$ be a radial weight. Then the following statements are true.
\begin{enumerate}
\item For each $z \in \bbD$, the function $B_z^\omega(\xi)$ is analytic in a larger disk than $\bbD$;
\item The kernel at zero, $B_0^\omega$, is a non-zero constant function;
\item There exists $r>0$ such that kernels $B_z^\omega$ are zero-free, if $|z|<r$;
\item If $B_z^\omega$ has a zero and $|w|>|z|$, then the kernel $B_w^\omega$ also has a zero;
\item If $B_z^\omega$ has a zero for some $z$, then the function $B^\omega(z,\xi)=B_z^\omega(\xi)$ vanishes on infinite number of points $(z,\xi)$ in the bi-disk $\bbD \times \bbD$;
\item No function $B_z^\omega$ can have infinitely many zeroes in the disk.
\end{enumerate}
\ENT

\begin{proof}
Let $z \in \bbD$ and $|z|<|w|<1$. By using \eqref{srs}, we may write
$$B_z^\omega(\xi)=\sum_{n=0}^\infty \frac{1}{\omega_n}(\overline{z}\xi)^n=\sum_{n=0}^\infty \frac{1}{\omega_n}(\overline{w(z/w)}\xi)^n=B_w^\omega(\overline{z/w}\xi).$$
The rightmost series represents an analytic function on $\{|\xi|<|w/z|\}$. This proves $(1)$.

The condition $(2)$ follows directly from the series representation of $B_z^\omega$. Since $\omega$ is assumed to be integrable, the number $\omega_0$ is finite.

Write
$$B_z^\omega(\xi)=\frac{1}{\omega_0}+\overline{z}\xi\sum_{n=1}^\infty \frac{1}{\omega_n}(\overline{z}\xi)^{n-1}.$$
By $(1)$ the function represented by the sum in the right-hand-side is analytic on a larger disk than $\bbD$, in particular, bounded. By choosing $|z|$ small enough, we can guarantee that
$$\frac{1}{\omega_0}>\left|\overline{z}\xi\sum_{n=1}^\infty \frac{1}{\omega_n}(\overline{z}\xi)^{n-1}\right|,$$ when $\xi \in \bbD$,
which implies $(3)$.

If $B_z^\omega(\xi_0)=0$, then $B_w^\omega(\overline{(z/w)}\xi_0)=0$, which proves $(4)$.

If $B_z^\omega(\xi_0)=0$, then there are infinitely many pairs $(w,\overline{(z/w)}\xi_0) \in \bbD\times \bbD$ (at least one for each $|w|>|z|$), for which $B^\omega$ vanishes, so $(5)$ holds.

Finally, suppose that $B_z^\omega$ has infinitely many zeroes on $\bbD$. If $|z|<|w|<1$, then from $(4)$ it follows that $B_w^\omega$ has infinitely many zeroes in $\{|\xi|<|z/w|\}$. This would force $B_w^\omega$ to collapse to the zero function. This is impossible, as it would mean that all functions in $A^2_\omega$ would vanish at $w \in \bbD$, while on the other hand, we know that $A^2_\omega$ contains all constants if $\omega$ is integrable. We obtain $(6)$.
\end{proof}

\section{Comments on the Segal-Bargmann type spaces.}

We end this paper with a short discussion on how the results adapt to the Segal-Bargmann (some interesting discussion on the good choices for alternative names for such spaces can be found in \cite{JPR} page 63) spaces of the complex plane. Let $\gamma>0$ and $0<p<\infty$. The Segal-Bargmann space $F^p_\gamma$ consists of entire functions $f$ satisfying
$$\|f\|_{F^p_\gamma}:=\left(\gamma p/2 \int_\bbC |f(z)|^p e^{-(\gamma p/2)|z|^2}dA(z)\right)^{1/p}<\infty.$$ The norm on $F^\infty_\gamma$ is, as one would expect, given by
$$\|f\|_{F^\infty_\gamma}:=\sup_{z \in \bbC} |f(z)|e^{-(\gamma/2) |z|^2}.$$
The space $F^2_\gamma$ is a reproducing kernel Hilbert space with a reproducing kernel $B_z^\gamma(\xi)=e^{\gamma \overline{z}\xi}$. For technical reasons, the exponent $p$ is included in the weight, see \cite{JPR}. For instance, this is known to guarantee that the orthogonal projection arising from $p=2$ is bounded for all $1\leq p \leq \infty$. To adapt Theorem \ref{sharp} to the Segal-Bargmann case, note that for $0<p<\infty$ the $L^p$-norm of $F^p_\gamma$ actually comes from $F^2_{\gamma p/2}$. Thus, since the case $p=\infty$ is obvious, we can easily obtain the sharp bound
$$|f(z)|\leq e^{(\gamma/2)|z|^2}\|f\|_{F^p_\gamma},$$ which is the formula $(1.6)$ of \cite{GK}, see also Theorem 2.8 \cite{ZF}. To be completely honest, the case $0<p<\infty$ also needs part of the argument from \cite{Vuk}.

For the remaining theorems, we write $\omega_\gamma(z)=\gamma e^{-\gamma|z|^2}$, and

$$\omega_\gamma^*(z)=\int_{|z|}^\infty \omega_\gamma(s)\log\left(\frac{s}{|z|}\right)s ds$$ for its associated weight defined on $\bbC\setminus \{0\}$. We see, like in Theorem \ref{example}, that the corresponding reproducing kernel for would be a constant multiple of $(1+\gamma \overline{z}\xi)e^{\gamma \overline{z}\xi}$, which clearly has a zero. By iterating the process, as in the proof of Theorem \ref{prescribe}, we can obtain kernels with any finite number of zeroes by a simple application of the fundamental theorem of algebra. The proof is actually much simpler than in the case of the unit disk, where we needed to use Rouche's theorem instead. We omit the details.

Finally, if $W$ is a radial weight on the complex plane such that polynomials are square-integrable with respect to it, conditions analogous to $(2)$, $(4)$ and $(5)$ of Theorem \ref{inf} can be verified. However, the condition $(1)$ does not really make sense, while $(3)$ and $(6)$ obviously do not follow. In fact we do not know, whether it is possible to find a radial $W$ on $\bbC$ with $B_z^W$ having infinitely many zeroes.

\end{document}